\documentclass[preprint,12pt]{elsarticle}
\usepackage{amssymb,amsmath}
\usepackage[usenames,dvipsnames,svgnames,table]{xcolor}
\usepackage{color}
\newtheorem{thm}{Theorem}[section]
\newtheorem{cor}[thm]{Corollary}
\newtheorem{lem}[thm]{Lemma}

\newtheorem{defn}[thm]{Definition}

\newcommand{\R}{\mathbb{R}}
\newcommand{\N}{\mathbb{N}}
\newcommand{\Q}{\mathbb{Q}}
\newcommand{\Z}{\mathbb{Z}}
\newcommand{\C}{\mathbb{C}}

\newcommand{\X}{\mathbb{X}}
\newcommand{\Y}{\mathbb{Y}}
\usepackage{hyperref}

\journal{be published}

\begin{document}

\begin{frontmatter}



\title{Discontinuous Almost
Automorphic Functions and Almost Automorphic Solutions of Differential Equations with Piecewise Constant Argument}


\author[ACH]{Alan Ch\'avez}
\author[SC]{Samuel Castillo}
\author[ACH]{Manuel Pinto}
\address[ACH]{
Universidad de Chile.
Departamento de Matem\'aticas.
Facultad de Ciencias.\\
Casilla 653.
Santiago-Chile.
{\tt alancallayuc@gmail.com, pintoj.uchile@gmail.com}}
\address[SC]{
Universidad del B\'{\i}o-B\'{\i}o.
Departamento de Matem\'atica.
Facultad de Ciencias.\\
Casilla 5-C.
Concepci\'on-Chile.
{\tt scastill@ubiobio.cl}
}

\begin{abstract}
In this article we introduce a class of discontinuous almost
automorphic functions which appears naturally in the study of almost
automorphic solutions of differential equations with piecewise
constant argument. Their fundamental  properties are used to prove
the almost automorphicity of bounded solutions of a system of
differential equations with piecewise constant argument. Due to the
strong discrete character of these equations, the existence of a
unique discrete almost automorphic solution of a non-autonomous
almost automorphic difference system is obtained, for which
conditions of exponential dichotomy and discrete Bi-almost
automorphicity are fundamental.

\end{abstract}

\begin{keyword}
Almost automorphic functions \sep Differential equation with
piecewise constant argument \sep Difference equations \sep
Exponential dichotomy.


\end{keyword}

\end{frontmatter}


\section{Introduction.}
A first order differential equation with piecewise constant argument (DEPCA) is an equation of the type
$$x'(t)=g(t,x(t),x([t])),$$
where $[\cdot]$ is the greatest integer function. The study of DEPCA
began in the year 1983 with the works of S. M. Shah and J. Wiener
\cite{SHAW}, then in 1984 K. L. Cooke and J. Wiener studied DEPCA
with delay \cite{cooke1}. DEPCA are  of considerable importance in
applications to some biomedical dynamics, physical phenomena (see
\cite{MU3,SBU} and some references therein), discretization problems
\cite{JW1,ZHYXXW}, etc; consequently they have had a huge
 development, \cite{KSCP1,KSCP2,WD,GP1,GP2,Mpinto,HRPR} (and some references therein) are evidence of these fact. In this
way many results about existence, uniqueness, boundedness,
periodicity, almost periodicity, pseudo almost periodicity,
stability  and other properties of the solutions for these equations
have been developed (see \cite{LC,MU3,29,GP1,GP2,RYJL,CHZ} and some
references therein). Recently, in  2006 the study of the almost
automorphicity of the solution for a DEPCA
was considered \cite{WD,NT}.\\

Let $\X, \Y$ be Banach spaces and $BC(\Y;\X)$ denotes the space of
continuous and bounded functions from $\Y$ to $\X$. A function $f\in
BC(\R;\X)$ is said to be almost automorphic (in the sense of S.
Bochner) if given any sequence $\{s_n'\}_{n\in\N}$ of real numbers,
there exist a subsequence
$\{s_n\}_{n\in\N}\subseteq\{s_n'\}_{n\in\N}$ and a function $\tilde
f$, such that the pointwise limits
\begin{equation}\label{eqx1}
\lim_{n\to\infty}f(t+s_n)=\tilde f(t),\ \lim_{n\to\infty}\tilde f(t-s_n)=f(t), t\in\R
\end{equation}
hold.

When the previous limits are uniform in all the real line, we say
that the function $f$ is almost periodic in the Bochner sense.
Following the classical notation we denote by $AP(\R;\X)$ and
$AA(\R;\X)$ the Banach spaces  of almost periodic and almost
automorphic functions respectively. For detailed information about
these functions we remit to the references
\cite{5aaa,5aa,5bb,CC,29,29GM,NVMET}. \\

\noindent Our interest in this work is to prove the almost
automorphicity of the bounded solutions of the DEPCA
\begin{equation}\label{eq1}
x'(t)=Ax(t)+Bx([t])+f(t),
\end{equation}

\noindent where $A,B \in M_{p\times p}(\R)$ are matrices and $f$ is an almost automorphic function. \\

The following definition expresses what we understand by solution for the DEPCA  (\ref{eq1}).

\begin{defn}\label{SOLDEPCA}{\rm
A function $x(t)$ is a solution of a DEPCA (\ref{eq1}) in the interval $I$, if the following conditions are satisfied:
\begin{itemize}
\item [i)]$x(t)$ is continuous in all $I$.
\item [ii)]$x(t)$ is differentiable in all $I$, except possibly in the points $n\in I\cap\Z$ where there should be a
lateral derivative.
\item [iii)]$x(t)$ satisfies the equation in all the open interval  $]n,n+1[, n\in \Z$ as well as is satisfied by its right
side derivative in each $n\in \Z$.
\end{itemize}}
\end{defn}

DEPCA are differential equations of hybrid type, that is, they have the structure of continuous and discrete dynamical systems, more precisely in (\ref{eq1}) the
continuity occurs on intervals of the form  $]n,n+1[, n\in\Z$ and the discrete aspect on $\Z$. Due to the continuity of the solution on the whole line for a  DEPCA,
we get a recursion formula in $\Z$ and thus, we can pass from an interval to its consecutive. The recursion formula appears naturally as solution of a
difference equation. \\

With this objective, we study a general non-autonomous difference equation
\begin{equation}\label{eq2}
x(n+1)=D(n)x(n)+h(n), n\in \Z,
\end{equation}
where $D(n)\in M_{p\times p}$ is a discrete almost automorphic
matrix and $h$ is a discrete almost automorphic function. To study
the equation (\ref{eq2}) we use conditions of exponential dichotomy
and a Bi-almost automorphic Green function \cite{Mpinto001,TJXX},
obtaining a theorem about the existence of a unique discrete almost
automorphic solution for (\ref{eq2}). In \cite{Mpinto001,TJXX},
functions with a Bi-property have shown to be very useful. When
$D(n)$ is a constant operator on an abstract Banach space, D. Araya
et al. \cite{38} obtained the existence of discrete almost
automorphic solutions under some geometric assumptions on the Banach
space and spectral conditions on the
operator $D$.\\

Note that although an almost automorphic solution $x$ of (\ref{eq1})
is continuous, the function $x([t])$  does not and then it is not
almost automorphic. Really, $x([t])$ has friendly properties for our
study when in (\ref{eqx1}), $\{s_n\}_{n\in\N}$ are in $\Z$. This
class of discontinuous functions, which we call $\Z$-almost
automorphic (see definition 2.1), appears inevitably in DEPCA and
allow us to study almost automorphic DEPCAs in a correct form (see
the notes about Theorem \ref{TEQ1}). This type of problem is present
in the study of continuous solutions of DEPCA of diverse kind as
periodic or almost periodic type, but it is not sufficiently
mentioned in the literature (see
\cite{MU2,MU3,Nieto2,Nieto3,Nieto1,RY,RYRZ}). The treatment of
almost periodic solutions for a DEPCA was initiated by R. Yuan and
H. Jialin \cite{RYJL}. E. A. Dads and L. Lachimi \cite{LC}
introduced discontinuous almost periodic functions to study the
existence of a unique pseudo almost periodic solution in a well
posed form to a DEPCA with delay. $\Z$-almost automorphic functions
generalize the ones proposed in \cite{LC}.\\

Properties derived in Section 2 for $\Z$-almost automorphic
functions allow us to simplify the proofs of some important results,
some of them known for almost automorphic functions in the
literature (see Theorem \ref{TD3} and \cite[Lemma 3.3]{NT}). We will
see that to obtain almost automorphic solutions of DEPCAs is
sufficient to consider $\Z$-almost automorphic perturbations. An
application of these facts is given by the use of the reduction
method in DEPCA (\ref{eq1}). This paper is organized as follows. In
Section 2, we introduce the $\Z$-almost automorphic functions and
their basic properties. In Section 3, we introduce the discrete
Bi-almost automorphic condition for the Green matrix to study
discrete non-autonomous almost automorphic solutions. Finally, in
Section 4, we study the almost automorphic solutions of equation
(\ref{eq1}) in several cases.

\section{$\Z$-almost automorphic functions.}

In this section we specify the definition of $\Z$-almost automorphic
functions with values in $\C^p$ and develop some of their
fundamental properties. Let us denote by $B(\R;\C^p)$ and
$BC(\R;\C^p)$ the Banach spaces of respectively bounded and
continuous bounded functions from $\R$ to $\C^p$ under the norm of
uniform convergence. Now  define $BPC(\R,\C^p)$ as the space of
functions in $B(\R;\C^p)$ which are continuous in $\R\backslash \Z$
with finite lateral limits in $\Z$. Note that $BC(\R;\C^p) \subseteq
BPC(\R,\C^p)$.

\begin{defn}
\label{2.1} A function $f\in BPC(\R;\C^p)$  is said to be
$\Z$-almost automorphic, if for any sequence of integer numbers
$\{s_n'\}_{n\in\N} \subseteq \Z$ there exist a subsequence
$\{s_n\}_{n\in\N}\subseteq \{s_n'\}_{n\in\N}$ such that the
pointwise limits in (\ref{eqx1})  hold.
\end{defn}

When the convergence in Definition \ref{2.1} is uniform, $f$ is
called $\Z$-almost periodic. We denote the sets of $\Z$-almost
automorphic (periodic) functions by $\Z AA(\R;\C^p)$ $(\Z
AP(\R;\C^p))$. $\Z AA(\R;\C^p)$ is an algebra over the field $\R$ or
$\C$ and we have respectively $AA(\R;\C^p)\subseteq \Z AA(\R;\C^p)$
and $AP(\R;\C^p)\subseteq \Z AP(\R;\C^p)$. Notice that a $\Z$-almost automorphic function is locally integrable.\\

\noindent For functions in $BC(\R\times \Y;\X)$ we adopt the following notion of almost automorphicity.
\begin{defn}
A function $f\in BC(\R\times\Y;\X)$ is said to be almost automorphic uniformly in compact subsets of $\Y$, if given any compact set $K \subseteq \Y$ and a sequence $\{s_n'\}_{n\in\N}$ of real numbers, there exists a subsequence $\{s_n\}_{n\in\N}\subseteq\{s_n'\}_{n\in\N}$ and a function $\tilde f$, such that for all $\ x\in K$ and each $t\in\R$ the limits
\begin{equation}
\label{EX4}
\lim_{n\to\infty}f(t+s_n,x)=\tilde f(t,x),\ \lim_{n\to\infty}\tilde f(t-s_n,x)=f(t,x),
\end{equation}
hold.
\end{defn}
The vectorial space of almost automorphic functions uniformly in
compact subsets is denoted by $ AA(\R\times\Y;\X)$,
see \cite{29,29GM}.

\begin{lem}\label{LEM1}If $f\in AA(\R;\C^p)$ $(resp. \ \, AP(\R;\C^p))$, then $f([\cdot])\in \Z AA(\R;\C^p)$ $(resp.\ \, \Z AP(\R;\C^p)).$
\end{lem}

All the next results for $\Z AA(\R;\C^p)$ are also valid for $\Z AP(\R;\C^p)$.


\begin{lem}\label{LEM8} The space $\Z AA(\R;\C^p)$ is a Banach space under the norm of uniform convergence.
\end{lem}
\paragraph{Proof} We only need to prove that the space $\Z AA(\R;\C^p)$ is closed in the space of bounded functions under the topology of uniform convergence.
Let $\{f_n\}_{n\in\N}$ be a uniformly convergent sequence of
$\Z$-almost automorphic functions with limit $f$. By definition each
function of the sequence is bounded and  piecewise continuous with
the same points of discontinuities, it is not difficult to see that
the limit function $f$ is bounded and piecewise continuous. Given a
sequence $\{s_n'\}_{n\in\N} \subseteq \Z$, it only rest to prove the
existence of a subsequence $\{s_n\}_{n\in\N} \subseteq
\{s_n'\}_{n\in\N}$ and a function $\tilde f$, where the pointwise
convergence given in (\ref{eqx1}) holds. As in the standard  case of
the almost automorphic functions the approach follows across the
diagonal procedure,
see \cite{29,29GM}. $\square$\\

\begin{lem}\label{LEM9}
 Let $G:\C^p\to \C^p$ be a continuous function and $f\in \Z AA(\R; \C^p)$, then  $G(f(\cdot))\in \Z AA(\R; \C^p)$.
\end{lem}

\begin{lem}\label{LEM3} Let $f\in AA(\R\times\C^p;\C^p)$ and uniformly continuous on compact subsets of $\C^p$, $\psi\in \Z AA(\R;\C^p)$.  Then $f(\cdot,\psi(\cdot))\in \Z AA(\R;\C^p)$.
\end{lem}
\paragraph{Proof} We have that the range of $\psi\in\Z AA(\R;\C^p)$ is relatively compact, that is, $K=\overline{\{\psi(t),t\in \R\}}$ is compact.
Let $\{s_n'\}_{n\in\N}\subseteq\Z$ be an arbitrary sequence, then
there exist a subsequence $\{s_n\}_{n\in\N}\subseteq
\{s_n'\}_{n\in\N}$ and functions $\tilde{f}$ and $\tilde{\psi}$ such
that the  pointwise limits in (\ref{EX4})
and \[ \lim_{n\to+\infty}\psi(t+s_n)=\tilde \psi(t),
\lim_{n\to+\infty}\tilde\psi(t-s_n)=\psi(t), t \in \R
\] hold. The equality $\displaystyle\lim_{n\to+\infty}f(t+s_n,\psi(t+s_n))=\tilde
f(t,\tilde\psi(t))$ follows from
\begin{eqnarray*}
|f(t+s_n,\psi(t+s_n))-\tilde f(t,\tilde\psi(t))|&\leq&|f(t+s_n,\psi(t+s_n))-f(t+s_n,\tilde\psi(t))|+\\
&+&|f(t+s_n,\tilde\psi(t))-\tilde f(t,\tilde\psi(t))|.
\end{eqnarray*}
The proof of $\displaystyle\lim_{n\to+\infty}\tilde
f(t-s_n,\tilde\psi(t-s_n))= f(t,\psi(t))$ is analogous.
$\square$\\

\noindent With analogous arguments we can prove the following Lemma.
\begin{lem}\label{LEM4}Let $f\in AA(\R\times\C^p\times\C^p;\C^p)$ be uniformly continuous on compact subsets of $\C^p\times\C^p$,
$\psi \in AA(\R;\C^p)$, then $f(\cdot,\psi(\cdot),\psi([\cdot])) \in \Z AA(\R;\C^p).$
\end{lem}

\noindent Now we want to give a necessary condition to say when a $\Z$-almost automorphic function is almost automorphic.

\begin{lem}\label{LEM6}
 Let $f$ be a continuous $\Z$-almost automorphic (periodic) function. If $f$ is uniformly continuous in $\R$, then $f$ is almost automorphic (periodic).
\end{lem}
\paragraph{Proof} Let $\{s_n'\}_{n\in\N}$ be an arbitrary sequence of real numbers, then  there exists a subsequence $\{s_n\}_{n\in\N} \subseteq \{s_n'\}_{n\in\N}$ of the form $s_n=t_n+\xi_n$ with $\xi_n \in \Z$ and $t_n \in [0,1[$ such that $\displaystyle\lim_{n\to \infty}t_n=t_0 \in[0,1]$. Moreover, $\{\xi_n\}_{n\in\N}$ can be chosen such that the pointwise limits
\begin{equation}\label{eqx2}
\displaystyle\lim_{n\to\infty}f(t+\xi_n)=: g(t),\ \displaystyle\lim_{n\to\infty} g(t-\xi_n)=f(t),\; t\in \R
\end{equation}
hold. As $f$ is uniformly continuous, the function $g$ is too. Let
us consider
\begin{eqnarray*}
 |f(t+t_n+\xi_n)- g(t+t_0)|&\leq& |f(t+t_n+\xi_n)-f(t+t_0+\xi_n)|+\\
&+&|f(t+t_0+\xi_n)-g(t+t_0)|.\nonumber
\end{eqnarray*}
Let $\epsilon >0$,  $\delta=\delta(\epsilon)$ be the parameter in
the uniform continuity of $f$. Let $N_0=N_0(\epsilon)\in\N$ be such
that for every $n\ge N_0$, $|t_n-t_0|<\delta$. Then the uniform
continuity of $f$ ensures that
$|f(t+t_n+\xi_n)-f(t+t_0+\xi_n)|<\frac{\epsilon}{2}$. Moreover, by
(\ref{eqx2}) there exists $N_0'=N_0'(t,\epsilon)$ such that if
$n\geq N_0',$ then $ \ |f(t+t_0+\xi_n)- g(t+t_0)|<\frac{\epsilon}{2}
$. Therefore, given $n\geq M_0=\max \{N_0,N_0'\}$, we have
$$|f(t+s_n)-g(t+t_0)|<\epsilon.$$
Similarly, from
the uniform continuity of $g$ and (\ref{eqx2}) we conclude
that  $|g(t+t_0-s_n)-f(t)|<\epsilon$, for all $n\geq M_0$.
Then $f \in AA(\mathbb{R},\mathbb{C}^p)$. $\square$

\begin{lem}\label{LEM7}
 Let $f\in \Z AA(\R;\C^p)$ $(resp.$ $\Z AP(\R;\C^p)$. The function $F(t)=\int_0^t f(s)ds$ is bounded if and only if  $F(\cdot)$ is almost automorphic $(resp.$ almost periodic$)$.
\end{lem}
\paragraph{Proof} The proof of the sufficient condition is immediate. For the necessary condition, since $F$ is uniformly continuous, we need to prove that $F$ is $\Z$-almost automorphic, which follows by the same arguments of \cite[Theorem 2.4.4]{29}. $\square$\\

\begin{lem}\label{THF1}
Let $\Phi:\R \to \mathbb M_{p\times p}(\R) $ be an absolutely
integrable matrix and $A\in M_{p\times p}(\R)$ be a constant matrix.
The operators
$$(L f)(t)=\int_{-\infty}^{\infty}\Phi(t-s)f(s)ds\ {\rm and } \ (\Upsilon f)(t)=\int_{[t]}^{t}e^{A(t-s)}f(s)ds,$$
map $\Z AA(\R;\C^p)$ into itself.
\end{lem}
\paragraph{Proof} We only prove the Lemma for $L$, the proof for $\Upsilon$ is analogous. It is easy to see that the operator $L$ is bounded. Let $\{ s_n'\}_{n\in\N}$ be a sequence of integers. Since $f\in \Z AA(\R;\C^p)$, there exists a subsequence $\{ s_n\}_{n\in\N}\subseteq \{ s_n'\}_{n\in\N}$ and a function $\tilde f$ such that we have the pointwise limits in (\ref{eqx1}).
Define the function $g(t)=(L\tilde{f})(t)$. Then, by the Lebesgue
Convergence Theorem
$$\lim_{n \to +\infty}(Lf)(t+s_n)=\lim_{n \to +\infty}\int_{-\infty}^{\infty}\Phi(t-s)f(s+s_n)ds=g(t).$$
Analogously, the limit $\displaystyle\lim_{n\to\infty} g(t-s_n)=(Lf)(t)$ holds. $\square$\\

\section{Almost Automorphic Solutions of Difference Equations.}

As it is noted in the literature \cite{LC,NT,RYJL,CHZ}, difference equations are very important in DEPCA studies. In
this section, we are interested in obtaining discrete almost
automorphic solutions of the system
\begin{equation}\label{ED2}
x(n+1)=C(n)x(n)+f(n), \ n \in \Z,
\end{equation}
where 
$C(\cdot)\in M_{p\times p}(\R)$ is a discrete almost automorphic matrix and  $f(\cdot)$ is a discrete almost automorphic function.

\begin{defn}Let $\X$ be a Banach space. A function $f:\Z \to \X$ is called discrete almost automorphic, if for any sequence $\{s_n^{'}\}_{n\in\N} \subseteq \Z$, there exists a subsequence  $\{s_n\}_{n\in\N} \subseteq\{s_n^{'}\}_{n\in\N}$, such that the following pointwise limits
$$\lim_{n \to +\infty}f(k+s_n)=:\tilde f(k), \lim_{n \to +\infty}\tilde f(k-s_n)= f(k), \ k \in \Z$$
hold.
\end{defn}
We denote the vector space of almost automorphic sequences by
$AA(\Z,\X)$ which becomes a Banach algebra over $\R$ or $\C$ with
the norm of uniform convergence (see \cite{38}). In
\cite{Mpinto001,TJXX}, we see the huge importance of the Bi-property
of a function $H:=H(\cdot,\cdot)$, such as Bi-periodicity, Bi-almost
periodicity, Bi-almost automorphicity; i.e. $H$ has simultaneously
the property in both variables. This motives the following
definition.

\begin{defn} For $\X$ being a Banach space, a function $H:\Z\times \Z \to \X$ is said to be a discrete Bi-almost  automorphic function, if for any sequence $\{s_n^{'}\}_{n\in\N} \subseteq \Z$, there exists a subsequence $\{s_n\}_{n\in\N} \subseteq\{s_n^{'}\}_{n\in\N}$, such that the following pointwise limits
$$\lim_{n \to +\infty}H(k+s_n,m+s_n)=:\tilde H(k,m), \lim_{n \to +\infty}\tilde H(k-s_n,m-s_n)= H(k,m), \ k,m \in \Z$$
hold.
\end{defn}
Some examples of discrete Bi-almost automorphic functions can be obtained by restriction to the integer numbers of continuous Bi-almost automorphic (periodic) functions in $\R$.

The following definition deals with the discrete version of exponential dichotomy \cite{CHZ}. Suppose that the matrix function $C(n),n\in \Z$, of the equation (\ref{ED2}) is invertible and consider $Y(n),n\in\Z$, a fundamental matrix solution of the system
\begin{equation}\label{DICHEDIS}
x(n+1)=C(n)x(n),n\in\Z.
\end{equation}

\begin{defn}
The equation (\ref{DICHEDIS}) has an exponential dichotomy with parameters  $(\alpha,K,P)$, if there are  positive constants $\alpha,K$ and a projection $P$ such that
$$|G(m,l)|\leq Ke^{-\alpha|m-l|}, m,l \in \Z,$$
where $G(m,l)$ is the discrete Green function which takes the explicit form
$$
G(m,l): = \left\{
\begin{array} {l}
Y(m)PY^{-1}(l),m\geq l\\
-Y(m)(I-P)Y^{-1}(l),m< l.
\end{array}
\right.
$$
\end{defn}

Now, we  give conditions to obtain a unique discrete almost automorphic solution of the system (\ref{ED2}).

\begin{thm}\label{TD1}Let $f\in AA(\Z,\C^p)$. Suppose that the homogeneous part of the equation (\ref{ED2})  has an
($\alpha$,K,P)-exponential dichotomy with discrete Bi-almost automorphic Green function $G(\cdot,\cdot)$. Then the
unique almost automorphic solution of (\ref{ED2}) takes the form:
\begin{equation}\label{eqP1}
 x(n)=\sum_{k \in \Z}G(n,k+1)f(k), \ n \in \Z
\end{equation}
and
$$
 |x(n)|\leq K(1+e^{-\alpha})(1-e^{-\alpha})^{-1}||f||_{\infty}, \, n \in \Z.
$$
\end{thm}

\paragraph{Proof} It is well known that the function given by (\ref{eqP1}) is the unique bounded solution of the discrete equation (\ref{ED2}) (see \cite [Theorem 5.7]{CHZ}). We prove that this solution is discrete almost automorphic. In fact, consider an arbitrary sequence $\{s_n'\}_{n\in\N} \subseteq \Z$. Since $f\in AA(\Z,\C^p)$ and $G(\cdot,\cdot)$ is discrete Bi-almost automorphic, there are a subsequence $\{s_n\}_{n\in\N} \subseteq\{s_n'\}_{n\in\N}$ and
functions $\tilde f(\cdot), \tilde G(\cdot,\cdot)$ such that the following pointwise limits
$$\lim_{n\to +\infty}f(m+s_n)=:\tilde f(m),  \lim_{n\to +\infty}\tilde f(m-s_n)= f(m), m\in  \Z$$
and
$$\lim_{n\to +\infty}G(m+s_n,l+s_n)=:\tilde G(m,l), \lim_{n\to +\infty}\tilde G(m-s_n,l-s_n)= G(m,l), m,l\in
 \Z$$
hold. Note that $|\tilde G(m,l)|\leq Ke^{-\alpha |m-l|}, \ m,l \in \Z$. Then,
\begin{eqnarray*}
x(n+s_n)&=&\sum_{k \in \Z}G(n+s_n,k+1)f(k)\\
&=&\sum_{k \in \Z}G(n+s_n,k+1+s_n)f(k+s_n),
\end{eqnarray*}
and from the Lebesgue Dominated Convergence Theorem we conclude that
$$\lim_{n\to\infty}x(n+s_n)=\tilde x(n),$$
where
$$\tilde x(n)=\sum_{k \in \Z}\tilde G(n,k+1)\tilde f(k).$$
To demonstrate the limit
$$\lim_{n\to\infty}\tilde x(n-s_n)=\sum_{k \in \Z}G(n,k+1) f(k)=x(n),$$
we proceed analogously.  $\square$\\

\section{Almost Automorphic Solutions for Linear  DEPCA.}

\noindent Finally, in this section we investigate the almost automorphic solution of the equation (\ref{eq1}). Before that, we reproduce the following useful result.

\begin{lem}\label{AUXLEM1}
 Let $f(\cdot)$ be a locally integrable and bounded function. If $x(\cdot)$ is a bounded solution of (\ref{eq1}), then $x (\cdot)$ is uniformly continuous.
\end{lem}
\paragraph{Proof} Since $x(\cdot)$ and $f(\cdot)$ are bounded, there is a constant $M_0>0$ such that
$\displaystyle\sup_{u\in\R}|Ax(u)+Bx([u])+f(u)|\leq M_0$. Now, as a
consequence of the continuity of $x$, we conclude that
$$|x(t)-x(s)|\leq\left|\int_s^t(Ax(u)+Bx([u])+f(u))du\right|\leq M_0|t-s|.$$
Then, the Lemma holds. $\square$\\

For a better understanding, we study the equation (\ref{eq1}) in several cases.\\

\noindent {\bf 4.1.  B=0.} \\

In this case the equation (\ref{eq1}) becomes the system of differential equations
\begin{equation}\label{DEPCA2}
x'(t)=Ax(t)+f(t),
\end{equation}
which has been well studied when $f\in AA(\R;\C^p)$, see
\cite{29,NVMET}. But when $f\in\Z AA(\R;\C^p)$ we have the following
Massera type extension
\begin{thm}\label{THF2}
 Let $f\in \Z AA(\R;\C^p)$. If the eigenvalues of $A$ have non trivial real part, then the equation (\ref{DEPCA2}) has
a unique almost automorphic solution.
\end{thm}
\paragraph{Proof} Since the eigenvalues of $A$ have non trivial real part, it is well known that the system $x'(t)=Ax(t)$
has an exponential dichotomy, that is, there are projections $P,Q$
with $P+Q=I$ such that the bounded solution of (\ref{DEPCA2}) has
the form
\begin{eqnarray*}
x(t)
&=&\int_{-\infty}^{t}e^{A(t-s)}Pf(s)ds-\int_{t}^{+\infty}e^{A(t-s)}Qf(s)ds.\label{eqSNG5}
\nonumber
\end{eqnarray*}
By Lemma \ref{THF1} we can see that this solution is bounded and $\Z$-almost automorphic. By the following Lemma
\ref{LEM6}, we only need to show that this solution is uniformly continuous, but this is a consequence of Lemma
\ref{AUXLEM1}. The conclusion holds. $\square$\\

\noindent For the scalar equation
\begin{equation}\label{DEPCA1}
x'(t)=\alpha x(t)+f(t),
\end{equation}
Theorem \ref{THF2} implies the following result.
\begin{cor}
 Let $f\in \Z AA(\R;\C)$ and, the real part of $\alpha$, $\Re(\alpha) \not=0$. Then the scalar equation (\ref{DEPCA1})
has a unique almost automorphic solution, given by
\begin{eqnarray*}
 x_1(t)&=&\int_{-\infty}^te^{\alpha(t-s)}f(s)ds, \ {\rm for} \  \Re(\alpha)<0,\\
 x_2(t)&=&-\int_t^{+\infty}e^{\alpha(t-s)}f(s)ds, \ {\rm for} \  \Re(\alpha)>0.
\end{eqnarray*}
\end{cor}
\begin{thm}
\label{4.4}
 Let  $\alpha$ be a purely imaginary complex number and $f\in \Z AA(\R;\C)$. If $x(\cdot)$ is a bounded solution of
(\ref{DEPCA1}) then
$x(\cdot)$ is almost automorphic.
\end{thm}
\paragraph{Proof.} Let  $\alpha=\theta i$, with $\theta \in \R$, then the solution of (\ref{DEPCA1}) is
$$x(t)=e^{\theta t i}x(0)+\int_0^te^{\theta (t-s)i}f(s)ds,\  t \in \R. $$
Since $x(\cdot)$ is bounded, we have that $\int_0^te^{i\theta (t-s)}f(s)ds$ is bounded and, by Lemma \ref{LEM7}, is
almost automorphic. Therefore $x(\cdot)$ is almost automorphic. $\square$\\

\noindent {\bf 4.2. B $\not=$0.}\\

\noindent By the variation of parameters formula, the solution of DEPCA (\ref{eq1}), for $t\in [n,n+1[$ and $n\in\Z$, satisfies

\begin{eqnarray}\label{DEPCA1D0}\nonumber
x(t)
&=&Z(t,[t])x([t])+H(t),
\end{eqnarray}
where $\displaystyle
Z(t,\tau)=e^{A(t-\tau)}+\int_{\tau}^{t}e^{A(t-s)}Bds$ and
$\displaystyle H(t)=\int_{[t]}^{t}e^{A(t-s)}f(s)ds$.

  By continuity of the solution $x$, if $t \to (n+1)^{-}$ we
obtain the difference equation

\begin{equation}\label{DEPCA1D}
x(n+1)=C(n)x(n) +h(n), n \in \Z,
\end{equation}
where $C(n)=Z(n+1,n)$ and $\displaystyle h(n)=\lim_{t \to (n+1)^{-}} H(t)$. By Lemma \ref{THF1}, $Z(t,[t])$
and $H(t)$ are $\Z$-almost automorphic functions, hence $C(n)$ and
$H(n)$ are almost automorphic sequences.

For the existence of the solution $x=x(t)$ of DEPCA (\ref{eq1}) on
all of $\R$, we assume that the matrix $Z(t,\tau)$ is invertible for all
$t,\tau\in [n,n+1]$ and all $n \in \mathbb{Z}$, see \cite{MU3,Mpinto,JW1}. This hypothesis will be needed
in the rest of the section. For example, when $A$ and $B$ are
diagonal matrices, we have that
\[
\begin{array}{rcl}
Z(t,\tau)&=&\displaystyle e^{A(t-\tau)}+B\int_{\tau}^t e^{A(t-s)}ds\\
\\
&=&\displaystyle e^{A(t-\tau)}\left[I+B\int_{0}^{t-\tau}
e^{-Au}du\right]
\end{array}
\]
is invertible if and only if the next assumption holds.

Assume that the eigenvalues $\lambda_A$ of $A$ and $\lambda_B$ of
$B$ satisfy for $u \in [0,1]$
\begin{equation}
\label{eigen} \left\{\begin{array}{rcl} \displaystyle
\frac{\lambda_B}{\lambda_A}[1-e^{-u\lambda_A}] \neq
-1,&\mbox{if}&\lambda_A \neq 0,\\ \lambda_B u \neq
-1,&\mbox{if}&\lambda_A=0.
\end{array}\right.
\end{equation}

As Theorem \ref{TEQ1} below will show, the existence on all of
$\mathbb{R}$ of the solutions of (\ref{eq1}) also follows from
condition (\ref{eigen}) when matrices $A$ and $B$ are simultaneously
triangularizable.

\begin{thm}\label{TD3} Let $x$ be a bounded solution of (\ref{DEPCA2}) with $f\in \Z AA(\R;\C^p)$. Then $x$ is  almost automorphic if and only if $x(n)$ in (\ref{DEPCA1D}) is discrete almost automorphic.
\end{thm}
\paragraph{Proof} If $x$ is an almost automorphic solution then restricting it to $\Z$, $x(n)$ is discrete almost automorphic. For $f\in \Z AA(\R;\C^p)$ and $x(n)$ being an almost automorphic sequence, the function $x$ given by (\ref{DEPCA1D0}) is well defined. The proof of the almost automorphicity of $x$ will follow at once if we prove its $\Z$-almost  automorphicity,  by Lemma \ref{AUXLEM1}.

\noindent Let us take an arbitrary sequence $\{s_n'\}_{n\in\N} \subseteq \Z$.
Then there are a subsequence   $\{s_n\}_{n\in\N}\subseteq \{s_n^{'}\}_{n\in\N}$,
functions $\tilde f$ and $\nu$ such that  the  limits in
(\ref{eqx1}) and
\begin{equation*}\label{eqaux.2}
\lim_{n \to +\infty}x(k+s_n)=\nu(k), \lim_{n \to +\infty}\nu(k-s_n)=x(n), \  k \in \Z
\end{equation*}
hold. Now, consider the limit function
$$y(t)=Z(t,[t])\nu([t]) +\int_{[t]}^{t}e^{A(t-s)}\tilde f(s)ds.$$
Then,
$$|x(t+s_n)-y(t)|\leq|Z(t,[t])||x([t]+s_n)-\nu([t])|+\int_{[t]}^{t}|e^{A(t-s)}||f(t+s_n)-\tilde
f(s)|ds,$$ and for each $t \in \R$ we have  $\displaystyle\lim_{n\to
+\infty}x(t+s_n)=y(t)$. Analogously $\displaystyle\lim_{n \to
+\infty}y(t-s_n)=x(t).$ Then, the bounded solution $x$ is
$\Z$-almost automorphic.
$\square$\\

Note that, without using $\Z AA(\R;\C^p)$, to prove directly $x(n)\in AA(\Z;\C^p)$ implies $x \in AA(\R;\C^p)$ is much more difficult (see \cite[Lemma 3]{WD} and \cite[Lemma 3.3]{NT}).\\

\noindent {\bf 4.3. A = 0, B $\not=$0.}\\

Theorem \ref{TD3} includes this important case
\begin{equation}\label{DEPCA3}
x'(t)=Bx([t])+f(t),
\end{equation}
for which the existence condition is reduced to invertibility for
$t\in [0,1]$ of $I+tB$. Therefore the following result is obtained.

\begin{cor}\label{TD4}
Let $f\in \Z AA(\R;\C^p)$ and $x$ a bounded solution of (\ref{DEPCA3}). Then, $x$ is almost automorphic if and only if $x(n)$ is discrete almost automorphic.
\end{cor}

\noindent {\bf 4.4. Reduction Method.} \\

By ``simultaneous triangularizations" of matrices $A$ and $B$, we
understand that there is an invertible matrix, say $T$, which
simultaneously triangularizes  both matrices $A$ and $B$. There
exist various results to obtain conditions under which simultaneous
triangularization holds, see for example the monograph of Heydar
Radjavi and Peter Rosenthal \cite{HRPR} and some references therein.

\begin{thm}[Reduction Method]\label{TEQ1} Consider $f\in \Z AA(\R;\C^p)$ and suppose that the matrices $A,B$ of the system  (\ref{eq1})
have simultaneous triangularizations and satisfy (\ref{eigen}). Let
$x$ be a bounded solution of (\ref{eq1}), then $x$ is almost
automorphic if and only if $x(n),$ in (\ref{DEPCA1D}), is discrete
almost automorphic.
\end{thm}
\paragraph{Proof} If $x$ is almost automorphic, then its restriction to $\Z$ is discrete almost automorphic. We will
prove that if $x(n)$ is discrete almost automorphic, then $x(\cdot)$ is almost automorphic. In fact, since $A,B$ have a
simultaneous triangularization, there is an invertible matrix $T$ such that
$$T^{-1}AT=\bar A=\begin{bmatrix}
{\alpha_1}&{a_{12}}&{a_{13}}&{\cdots}&{a_{1p}}\\
{0}&{\alpha_2}&{a_{22}}&{\cdots}&{a_{2p}}\\
{\vdots}\\
{0}&{0}&{\cdots}&{0}&{\alpha_{p}}
\end{bmatrix}, T^{-1}BT=\bar B=\begin{bmatrix}
{\beta_1}&{b_{12}}&{b_{13}}&{\cdots}&{b_{1p}}\\
{0}&{\beta_2}&{b_{22}}&{\cdots}&{b_{2p}}\\
{\vdots}\\
{0}&{0}&{\cdots}&{0}&{\beta_{p}}
\end{bmatrix};$$
where, for $i\in \{1,2,\cdots,p\}$, $\alpha_i$ and $\beta_i$ are the
eigenvalues of $A$ and $B$ respectively. Consider the following
change of variables $y(t)=T^{-1}x(t),$ then the boundedness of
$x(t)$ is equivalent to the boundedness of $y(t)$, which is a
solution of the following new system
\begin{eqnarray*}
y'(t)&=&\bar A y(t)+\bar B y([t])+T^{-1}f(t).
\end{eqnarray*}
Observe that, by Lemma \ref{LEM9}, the sequence  $y(n)=T^{-1}x(n)
\in AA(\mathbb{Z},\mathbb{C}^p)$, since $x(n)$ is almost automorphic. Let
$T^{-1}f(t)=H(t)=(h_1(t),h_2(t),\cdots,h_p(t))$, then we have the
almost automorphic system
\begin{equation}\label{SIST1}
y'(t)=\bar A y(t)+\bar B y([t])+H(t),
\end{equation}
namely
\begin{eqnarray*}
y_1'(t)&=&\alpha_1
y_1(t)+\sum_{i=2}^{p}a_{1i}y_i(t)+\beta_1y_1([t])+\sum_{i=2}^{p}b_{1i}y_i([t])+h_1(t)\\
y_2'(t)&=&    \alpha_2y_2(t)+\sum_{i=3}^{p}a_{2i}y_i(t)+\beta_2y_2([t])+\sum_{i=3}^{p}b_{2i}y_i([t])+h_2(t)\\
&\vdots&\\
y_{p-1}'(t)&=& \alpha_{p-1}y_{p-1}(t)+a_{p-1p}y_p(t)+\beta_{p-1}y_{p-1}([t])+b_{p-1p}y_p([t])+h_{p-1}(t)\\
y_p'(t)&=&\alpha_py_p(t)+\beta_p y_p([t])+h_p(t).
\end{eqnarray*}
Now take the  $p$ th-equation
\begin{equation}\label{SIST12}
 y_p'(t)=\alpha_py_p(t)+\beta_p y_p([t])+h_p(t),
\end{equation}
where the eigenvalues $\alpha_p$ of $A$ and $\beta_p$ of $B$ satisfy
(\ref{eigen}).

Since $AA(\R;\C^p)\subseteq \Z AA(\R;\C^p)$ and $y_p$ is a bounded
solution of (\ref{SIST12}), from Theorem \ref{4.4},  $y_p(t)$ is
almost automorphic. Consider now the $(p-1)$ th-equation
$$y_{p-1}'(t)= \alpha_{p-1}y_{p-1}(t)+\beta_{p-1}y_{p-1}([t])+\left[a_{p-1p}y_p(t)+b_{p-1p}
y_p([t])+h_{p-1}(t)\right].$$

\noindent By Lemma \ref{LEM1}, $y_p([t])$ is $\Z$-almost automorphic, then the function
\[
z_{p-1}(t)=a_{p-1p}y_p(t)+b_{p-1p} y_p([t])+h_{p-1}(t)
\]
is again $\Z$-almost automorphic. Similarly, we can conclude that
$y_{p-1}(t)$ is an almost automorphic solution of the equation
\begin{equation}
\label{FEQU1} y_{p-1}'(t)=
\alpha_{p-1}y_{p-1}(t)+\beta_{p-1}y_{p-1}([t])+z_{p-1}(t),
\end{equation}
since it is a bounded solution. Following this procedure, we obtain the almost automorphic solution $y(t)$  of the system (\ref{SIST1}) and thus $x \in AA(\R,\C^p)$. $\square$\\

Note that the discontinuous function $z_{p-1}$ in (\ref{FEQU1}) is
$\Z$-almost automorphic, although functions  $h_{p-1}, h_p \in
AA(\R,\C)$. Then, the presence of $\Z$-almost automorphic terms is
proper of DEPCA. The $\Z$-almost automorphic space contains
correctly the $\Z$-almost periodic and the interesting $\Z$-periodic
situation (which are periodic functions not necessarily continuous),
see \cite{ACHMP}. Then we conclude.

\begin{cor}\label{corap}
Let $f\in \Z AP(\R,\C^p)$. Then, every bounded solution $x$ of the
DEPCA (\ref{eq1}) is almost periodic if and only if $x(n)\in
AP(\Z,\C^p)$.
\end{cor}

\begin{cor}\label{corperiodic} Suppose that $f$ is a $\Z$-$\omega$-periodic function, with $\omega \in\Q$, then\\
\noindent a) If $\omega=p_0\in \Z$, every bounded solution $x$ of the DEPCA (\ref{eq1}) is $\omega$-periodic if and only if the sequence $x(n), n\in\Z$, is discrete $\omega$-periodic. \\
\noindent b) If $\omega=\frac{p_0}{q_0}\in \Q$ with $p_0,q_0\in \Z$
relatively primes, then every bounded solution $x$ of the DEPCA
(\ref{eq1}) is $q_0\omega$-periodic if and only if the sequence
$x(n), n\in\Z$ is discrete $q_0\omega$-periodic.
\end{cor}

\section{Acknowledgments.}
This research has been partially supported by Fondecyt 1080034 and 1120709. The first author wants to thank
the kind hospitality of my friend Mr. Alvaro Corval\'an Azagra.
\bibliographystyle{plain}
\bibliography{biblio}
\nocite{*}


\end{document}